\documentclass[journal]{IEEEtran}
\IEEEoverridecommandlockouts  
\usepackage{cite}
\usepackage{caption}

\usepackage{subcaption}
\usepackage{stfloats}
\usepackage{amsmath,amssymb,amsfonts}
\usepackage{graphicx}
\usepackage{textcomp}
\usepackage{xcolor}
\usepackage{diagbox}
\usepackage{multirow}
\usepackage{makecell}
\usepackage{algorithm}
\usepackage{algpseudocode} 
\usepackage{amsmath}
\usepackage{tikz}
\usetikzlibrary{decorations.pathreplacing}
\usetikzlibrary{positioning}
\def\BibTeX{{\rm B\kern-.05em{\sc i\kern-.025em b}\kern-.08em
   T\kern-.1667em\lower.7ex\hbox{E}\kern-.125emX}}

\makeatletter
\IEEEtriggercmd{\reset@font\normalfont\fontsize{7.9pt}{8.40pt}\selectfont}
\makeatother
\IEEEtriggeratref{1}

\usepackage[capitalize]{cleveref} 
\begin{document}
\title{Existence of Trust-field in Vehicular Ad Hoc Networks: Empirical Evidence}
\author{
Md Mahmudul Islam, 
Shaurya Agarwal, \textit{Senior Member IEEE}

\thanks{Md Mahmudul Islam and Shaurya Agarwal are with the University of Central Florida, Orlando, FL, USA
{\tt\small shaurya.agarwal@ucf.edu}}}

\markboth{}%
{Islam \MakeLowercase{\textit{et al.}}: Existence of Trust-field in Vehicular Ad Hoc Networks: Empirical Evidence}
\maketitle
\begin{abstract}
Vehicular Ad Hoc Networks (VANETs) play a crucial role in enhancing road safety and traffic efficiency by enabling communication between vehicles (V2V) and between vehicles and infrastructure (V2I). Robust trust management is necessary to ensure the reliability of information in decentralized systems. This paper presents the notion of a ``Trust Field" in VANETs, conceptualized as the behavior of the nodes that represents trust levels evolving in both spatial and temporal dimensions. Using the LogitTrust model, we provide empirical evidence of how trust fields in vehicular networks change over time in different scenarios, including when malicious nodes are present. The results of our study demonstrate that the trust domain can adjust to fluctuations in network conditions, thereby offering a comprehensive metric for assessing the reliability of nodes. This innovative method improves the dependability of VANET applications by efficiently detecting and mitigating malicious actions.

\end{abstract}


\section{Introduction}
Vehicular ad-hoc networks (VANETs) are a special category of mobile ad-hoc networks (MANETs) that enable vehicles to communicate with each other in a decentralized fashion without the assistance of any infrastructure. VANETs were initially developed to improve road safety by enabling cooperative collision warning through vehicle-to-vehicle (V2V) and vehicle-to-infrastructure (V2I) communication. A vehicle can either report an emergency to other drivers on the road or seek information about the location of a gas station or a nearby parking lot \cite{hasrouny2017vanet}. VANETs are also utilized to alleviate the consequences of road collisions and provide prior warnings to drivers. Since the information communicated is disseminated in an environment where anybody can access it, ensuring security and privacy are paramount in VANETs \cite{tyagi2014investigating}. Therefore, every VANET must meet the security and privacy criteria to provide an effective and dependable system. It is mandatory to ensure the integrity of the exchanged messages to prevent any unauthorized insertion or modification by potential attackers, including insiders, outsiders, and malicious actors. Any impacted application might pose significant risks to drivers and passengers. According to the World Health Organization, road users' annual global fatality count exceeds 1.35 million \cite{who2018road}. Hence, authentication and trust are vital in ensuring secure and reliable data flow in VANETs.

In VANET, vehicles serve as both providers and consumers simultaneously. Therefore, effective data transmission among vehicles is crucial for successfully realizing VANET applications. The information and its source's credibility are essential to examine when deciding whether to accept or reject information obtained from neighbors. Before taking any action based on the received data, it is imperative to verify the quality and integrity of the information. Trust is a reliable means of dealing with vehicles that are performing maliciously \cite{hussain2020trust}. The dynamic and mobile nature of VANETs, combined with their open structure, renders them more vulnerable to internal and external attacks than traditional communication networks. The threats include denial of service (DoS), bogus information, and impersonation \cite{malhi2020security, mahmood2021security}. While cryptography provides hard security through data confidentiality and authentication, it falls short in detecting dynamic, malicious behaviors of nodes post-authentication. 

Trust management, therefore, becomes crucial for ensuring secure routing, efficient relay selection, and data reliability 
\cite{almarshoud2024security}. \textit{Trust} enables nodes to assess the reliability and accuracy of information from other nodes and complements cryptography by continuously evaluating node behavior and data genuineness. It leads to malicious node identification, credential revocation, and blackhole and DoS attack mitigation \cite{cui2019rsma, tolba2018trust, djamaludin2016revocation}. Furthermore, it can aid in privacy preservation and incentivize cooperation among vehicles, resulting in enhanced connectivity and network security 
\cite{wei2014security}.

A field is a physical quantity described by a scalar, vector, or tensor. It has a specific value assigned to every point in space and time through some functions. For instance, an electric field is an electrical property at every point in space where an electric charge is experiencing unit force by other charges. Similarly, a density field in a vehicular network is a vector field, where the density is a function of space and time, and at any point, it quantifies the concentration of vehicles. This paper defines trust in VANETs as a composite measure, integrating elements of data genuineness, agent reputation, and expectation of altruistic participation by agents. 

Although there has been much study on trust management in VANETs, the idea of a trust field that combines different trust measures and models to offer an extensive trust evaluation has not been thoroughly investigated. This research aims to address this shortcoming by examining the existence of a trust field in VANETs by providing some empirical findings.


%
%
%


\textbf{Contributions:} The primary contributions of this paper are as follows:
\textbf{(a) Conceptualization of the Trust Field:} We define and elaborate on the concept of a trust field in VANETs, integrating existing trust models and metrics into a cohesive framework.;
\textbf{(b) Empirical Analysis:} We present empirical evidence supporting the existence of a trust field through simulations and real-world data analysis, demonstrating its effectiveness in enhancing trust management in VANETs.
;

\textbf{Outline:} Section \ref{sec:related_works} provides a brief background on Trust computation and management and presents current methods for trust management in VANET. Section \ref{sec:trust_field} presents the system model, including dynamic trust and behavior pattern factors of the vehicular nodes, the attacker model, and the problem definition. Section \ref{sec:experiment} provides an overview of the LogitTrust model and how to evaluate trust using the model, followed by a discussion of the findings. Finally, section \ref{sec:conc} provides concluding remarks.
%
%
\section{State of the Art: Trust Computation and Management} \label{sec:related_works}
Modeling trust and its computation is challenging. In this section, we identify the state-of-the-art and research challenges.

\subsection{How is Trust Computed and Managed?}

Trust Management (TM) approaches quantify and predict trust, dynamically evaluating nodes and transmitted data.
They measure trust levels based on factors such as past behavior, shared information, and reputation. Trustworthiness evaluation can focus on the agents (entity-centric models) \cite{marmol2012trip, minhas2011multifaceted, soleymani2015trust}, the received data (data-centric models) \cite{hussain2016hybrid, rawat2015trust, gurung2013information}, or view trust through cooperation viewpoint \cite{shivshankar2015evolutionary, srinivasan2003cooperation}. Centralized TM approaches are susceptible to a single point of failure, carry high maintenance costs, and struggle with scalability \cite{li2012reputation, li2013rgte, bissmeyer2012central}. 
Decentralized methods, on the other hand, utilize information redundancy and consensus building  \cite{ltifi2016smart, sedjelmaci2015accurate, wahab2014cooperative, kumar2014collaborative}, requiring lightweight algorithms for timely decision-making. 
TM approaches include game-theoretic \cite{chen2016game, chiti2015context, wang2016game}, cryptography \cite{pham2018adaptive, hu2016replace, tangade2020trust}, and machine learning based \cite{oubabas2018secure, shams2018trust, fan2019trust}. Game theory, in particular, stands out given its capabilities in agent behavior analysis, clustering, incentive design, and convergence analysis. 

\subsection{Knowledge Gaps and Research Questions (RQs)} 

\noindent
\textbf{(RQ1)} \textit{\textbf{How to holistically capture trust perception in VANETs?}}

Existing trust computation methodologies gauge the degree to which an agent (or information) can be trusted under given conditions. However, they do not encapsulate \textit{trust perception} that an agent harbors towards the network based on other agents' behavior. An effective model should capture interactions among vehicles and be responsive to changes in network conditions and traffic operations. Efficient lifetime definitions for trust values are also essential due to storage constraints and the dynamic network topology. 

\vspace{2pt}

\noindent
\textbf{(RQ2)} \textit{\textbf{How to integrate traffic dynamics into trust management?}}
%
%


The effective modeling of trust evolution in V2V networks is \textit{{challenging}} due to (a) the dynamic and unpredictable network structure, which involves agents moving at various speeds, restrictions imposed by road geometry, and intermittent peer-to-peer communication, (b) the nonlinear interdependency and bidirectional relationship with the traffic flow. 
Understanding the spatial evolution of trust and its entanglement with the traffic flow can provide precious clues on how to effectively compute and manage trust in a timely and resource-aware fashion.
%
%
\section{Trust Computation Model} \label{sec:trust_field}
\subsection{Dynamic Trust}
We consider the concept of dynamic trust, which refers to the dynamic alteration of a vehicle's trust level as the operational and environmental conditions of the VANET change due to the mobility and behavior of the nodes. The trust level of each node will change dynamically based on current behaviors and past trust history. We have adapted and modified the Logit Trust method \cite{wang2014logittrust} to represent the trustworthiness of individual vehicles. To model the trust level of individual nodes, we considered two metrics derived from the VANET communication environment: (a) Packet Forwarding Ratio and (b) Packet Forwarding Delay.

\subsubsection{Packet Forwarding Ratio} This metric determines how many packets a node has forwarded to its neighboring vehicles compared to the number of packets it received. It helps identify a node's maliciousness if it intentionally drops packets to create black hole attacks \cite{PATIL2017317, kumar2021black}. For instance, if a vehicle receives $10$ packets and forwards $5$ packets, the packet forwarding ratio (PFR) is $0.5$, which can be defined as follows:

\begin{equation}
    PFR = \frac{P_f}{P_r}
\end{equation}

Where $P_f$ and $P_r$ represent the number of packets forwarded and received by a vehicle, respectively.

\subsubsection{Packet Forwarding Delay} This metric measures how long it takes a vehicular node to forward all the packets it receives during a sensing window. It helps find any vehicle that is trying to launch a packet delay attack \cite{prajapati2011implementation}. In this type of attack, a malicious node adds a delay to the forwarding of packets in the network. If a vehicle forwards \( P \) packets during a sensing window, the packet forwarding delay (PFD) can be calculated as follows:

\begin{equation}
    PFD = \sum_{p=1}^{P} \frac{1}{delay(p)}
\end{equation}

\subsection{Attaker Model}

In VANETs, malicious nodes attempt to disrupt normal operations by injecting falsified information or altering existing information. To model the behavior of such attackers, we designed the communication patterns so that they drop most of the packets to mimic black hole attacks while forwarding the remaining packets with some delay. Specifically, the malicious node drops $90\%$ of the packets it receives and forwards the rest with some delay ranging between $100$ and $500$ milliseconds.

\subsection{Problem Definition}

The problem is for vehicle $i$  to accurately determine if vehicle $j$ is transmitting accurate or fabricated data, based on a collection of past evidence. At time $t$, the observation $s_{ij}^t$ of vehicle $i$ represents the accuracy of the information received from vehicle $j$, which can be either true or false. If the information is accurate, vehicle $j$ is deemed reliable and represented by the value 1; otherwise, it is considered untrustworthy and represented by the value 0. The vehicle's behavior in VANET at time $t$ can be described by two communication metrics: PFR (Packet Forwarding Rate) and PFD (Packet Forwarding Delay). These metrics are represented by a column vector $\mathbf{x}^{t}$ = $[PFR, PFD]^\top$. Trust refers to the likelihood, denoted as $\theta_j^{t}$, that vehicle $j$ is transmitting accurate information to vehicle $i$ at time $t$, given the behavioral conditions indicated by $\mathbf{x}^{t}$.

Additionally, let $s_j$ represent the collection of evidence obtained by vehicle $i$ from its interaction with vehicle $j$ throughout the time interval $[0, T]$. Furthermore, let $\mathbf{x} = \{\mathbf{x}^t,t=1,...,T\}$ represent the behavior conditions over the interval $[0, T]$. LogitTrust addressed this issue by analyzing the behavioral pattern of vehicle $j$, represented by a latent variable $\beta_j$, between $s_j$ and $\mathbf{x}$. It then predicted $s_j^{T+1}$ based on $\mathbf{x}^{T+1}$, specifically using the expected value of $s_j^{t}$ given $\mathbf{x}^{t}$ and ${\beta}_j$, denoted as $ \mathbb{E}[s_j^{t} | \mathbf{x}^{t}, {\beta}_j]$. The conditional expectation is a real number that ranges from 0 to 1, indicating the trust level of the vehicle $j$ at time $T+1$ from vehicle $i$'s perspective.
\section{Hypothesis and Experimental Findings} \label{sec:experiment}

\subsection{Hypothesis: Existence of Trust Field}
It is critical to \textit{acknowledge} that the \textit{trust} evolves spatially as vehicles move through different locations, encountering various communication partners and environmental conditions that affect the perceived reliability and reputation of information sources. As vehicles traverse the network, they continuously update their trust assessments based on new interactions and the behavior of other nodes. For instance, a properly behaving trusted agent raises the trust score in its vicinity, indicating higher information reliability, while a malicious agent lowers the trust score, indicating potential unreliability. 
Recognizing that the evolution of trust in space and time is intertwined with vehicular flow dynamics, we \textit{\textbf{hypothesize}} that deciphering this relationship can serve as the \textit{linchpin} for understanding how trust disseminates across the road network as vehicles move. Moreover, supported by our preliminary results (see next section), we envisage that trust evolution can be conceptualized as a spatiotemporal field, hereinafter referred to as the ‘\textit{trust field}.’

\subsection{Trust Evaluation} \label{sec:trust_evaluation}
To evaluate the trust level, each vehicular node calculates the trust score of its neighbors at every delta timestep (ie. sensing window). Throughout the sensing window, vehicle $i$ collects the PFR and PFD for all of its neighbors and then calculate the trust score. For the very first interaction when there are no historical evidence available for vehicle $j$, vehicle $i$ uses equation $(\ref{eq:trust_first})$ to calculate the trust score $\mathbb{E}[s_j^{t} | \mathbf{x}^{t},{\beta}_j]$ and record the $\mathbf{x}^{t}$ and $s_j^{t}$ in the history table as evidence for later uses. 

\begin{equation}\label{eq:trust_first}
    \mathbb{E}[s_j^{t} | \mathbf{x}^{t}, {\beta}_j] = \theta_j^{t} = \frac{1}{1+{e^{-(\mathbf{x}^{T+1})^\top\beta_j}}}
\end{equation}

Where $\mathbf{x}^{t}$ is a column vector $[PFR, PFD]^\top$ representing the behavior of the vehicle $j$ in the VANET communication environment, $\theta_j^{t}$ is the expected trust score of the vehicle $j$ at time $t$, and $\beta_j$ is a column vector of coefficient. If $\theta_j^{t}$ is greater than the threshold value, then we consider the vehicle $j$ as trustworthy and denote $s_j^t = 1$; otherwise, $s_j^t = 0$.

From the following interaction, vehicle $i$ utilizes the previous history as evidence and current behavior of vehicle $j$ to calculate the trust level using Algorithm \ref{alg:logit_trust} and update at every delta timestep.

\begin{algorithm}
\caption{LogitTrust}
\label{alg:logit_trust}
\begin{algorithmic}[1]
\State \textbf{Input:} $\mathbf{x}, s_j, \nu_0, \mathbf{x}^{T+1}$
\State \textbf{Output:} $\mathbb{E}[s_j^{T+1} | \mathbf{x}^{T+1}, \hat{\beta}_j]$
\State $k \gets 0$
\State $\hat{\beta}_j^{(k)} \gets 1$
\While{not converged}
    \For{$t \gets 1$ to $T$}
        \State $u^t \gets (\mathbf{x}^t)^\top \hat{\beta}_j^{(k)}$
        \State $w^t \gets \frac{s_j^t - (2s_j^t - 1) F_{\nu_0 + 2} (-(1 + 2/\nu_0)^{1/2} u^t)}{s_j^t - (2s_j^t - 1) F_{\nu_0} (-u^t)}$
        \State $s_j^{t*} \gets u^t + \frac{(2s_j^t - 1) f_{\nu_0} (u^t)}{s_j^t - (2s_j^t - 1) F_{\nu_0 + 2} (-(1 + 2/\nu_0)^{1/2} u^t)}$
    \EndFor
    \State $S_0 \gets \sum_{t=1}^{T} w^t \mathbf{x}^t (\mathbf{x}^t)^\top$
    \State $S_1 \gets \sum_{t=1}^{T} w^t \mathbf{x}^t s_j^{t*}$
    \State $k \gets k + 1$
    \State $\hat{\beta}_j^{(k)} \gets S_0^{-1} S_1$
\EndWhile
\State \Return $\mathbb{E}[s_j^{T+1} | \mathbf{x}^{T+1}, \hat{\beta}_j] \gets \frac{1}{1 + \exp(- (\mathbf{x}^{T+1})^\top \hat{\beta}_j^{(k)})}$
\end{algorithmic}
\end{algorithm}

\subsection{Experimental setup}
\begin{figure}[ht]
	\centering
	\includegraphics[width=\columnwidth]{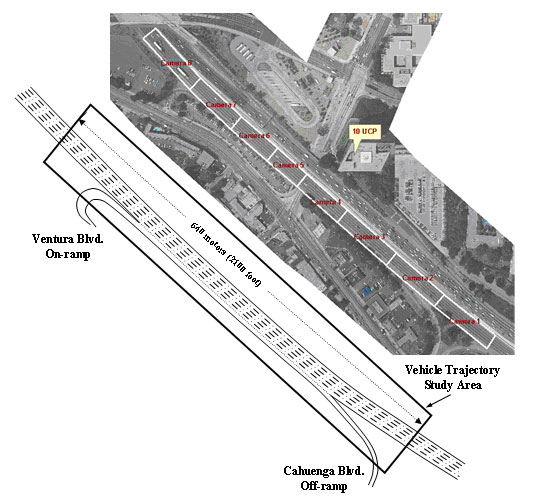}
	\caption{The schematic drawing shows the number of lanes and location of the on-ramp at Ventura Boulevard and the off-ramp at Cahuenga Boulevard within the US 101 study area. \cite{us_highway_101_dataset}}
	\label{fig:network_ngsim_us_101}
\end{figure}

In order to assess the performance of the LogitTrust model on actual vehicular traffic, we have selected the NGSIM US-101 highway vehicular trajectory data gathered by the US Federal Highway Administration. The NGSIM dataset offers detailed information on each vehicle's exact position, including its lane position and relative placement to other vehicles, recorded every one-tenth of a second on 2100 ft stretches of the southbound US-101 freeway. In essence, the NGSIM data offers an in-depth view of traffic by tracking individual vehicles rather than measuring the overall velocity or density of the macroscopic flow. The road network of the NGSIM US-101 is illustrated in Fig. \ref{fig:network_ngsim_us_101}.

\begin{figure}[ht]
	\centering
	\includegraphics[width=\columnwidth]{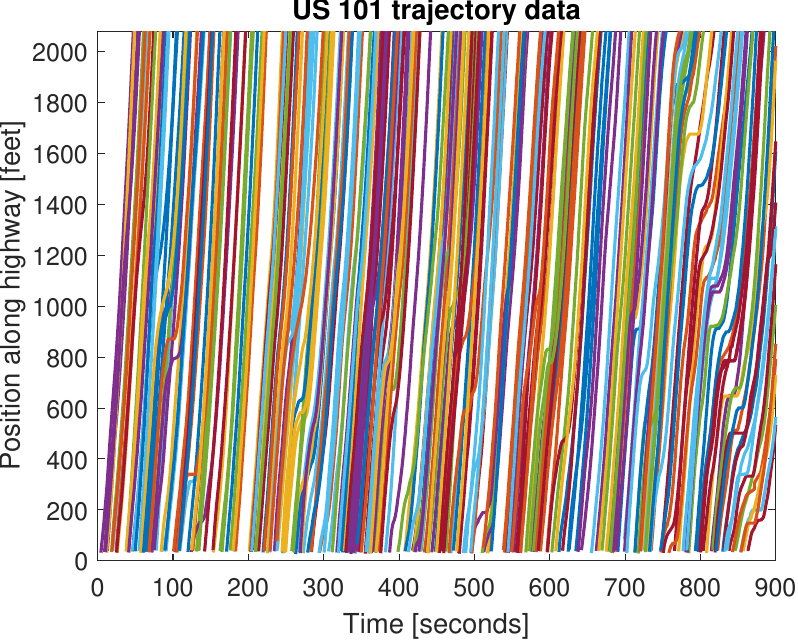}
	\caption{NGSIM US-101 trajectory data for the first 15 min. The data was collected between the hours of 7:50 am and 8:35 am, during the onset of congestion. The section of the highway studied consists of five main lanes, a single on and off-ramp and an auxiliary lane between the on and off-ramp. Every colored line corresponds to a unique vehicle.}
	\label{fig:trajectory_ngsim_us_101}
\end{figure}

The first 15 minutes of the vehicle trajectories for the NGSIM US-101 freeway are shown in Fig. \ref{fig:trajectory_ngsim_us_101}. Using the procedures mentioned earlier, we utilized the trajectories to calculate the trust score for each vehicular node at every delta timestep. We first converted these trajectories into the ns-3 mobility format and then imported into ns-3 for network simulation. In ns-3, we have designed a custom vehicular application where each vehicle transmits a message to its neighboring vehicles. Upon receiving the message, the nodes will forward it further to its neighbor, and the process keeps going until the message is forwarded to 3 hops. Then, we collect the network performance metrics: PFR and PFD for each vehicle. Finally, we implemented the Logit Trust algorithm to calculate the trust score of each vehicle at every delta timestep.

\begin{figure}[ht]
	\centering
	\includegraphics[width=\columnwidth]{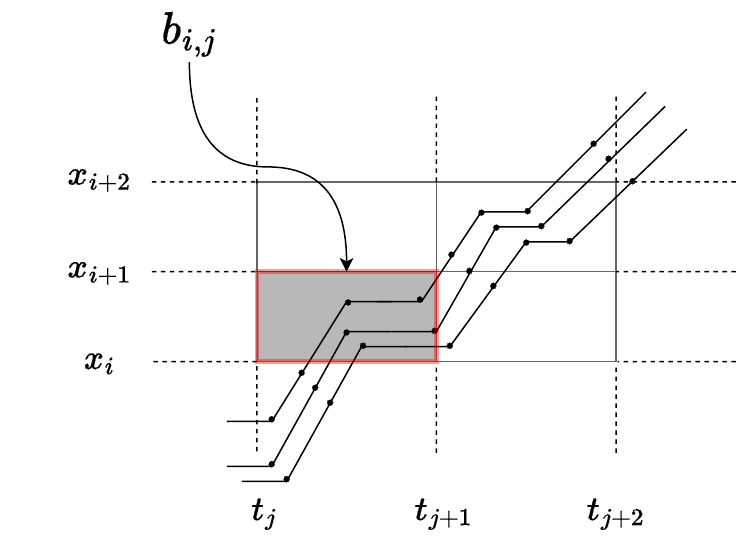}
	\caption{Schematic of the binning method utilized to generate the spatiotemporal data. The spatiotemporal domain is divided into bins in which the trust scores are computed according to the Eq. (\ref{eq:bin_trust}).}
	\label{fig:binning_method}
\end{figure}

The trust value we calculated is based on microscopic data associated with each vehicle. However, we are interested in extracting the trust field, which is a spatiotemporal pattern of the trust score given for every point in time and space. So, we are required to construct macroscopic data from the NGSIM trajectory. To accomplish this, we utilize the binning techniques established by \cite{avila2020data} and divide the spatiotemporal area $[0, L] \times [0, T]$ into discrete bins of dimensions $\Delta X \times \Delta T$, as depicted in Fig. \ref{fig:binning_method}. The variable $L$ denotes the overall length of the freeway, $T$ represents the total duration of data collection, $\Delta X$ indicates the size of each spatial step, and $\Delta T$ reflects the size of each temporal step. Hence, the formula for a single bin is expressed as follows in equation $(\ref{eq:binning})$, where $n_x = \frac{L}{\Delta X}$ and $n_t = \frac{T}{\Delta t}$ represent the number of bins in space and time, respectively.

\begin{equation}\label{eq:binning}
\text{Bin}_{i,j} = \left[ i \Delta X, (i + 1) \Delta X \right] \times \left[ j \Delta T, (j + 1) \Delta T \right]_{i \in (0, n_x), j \in (0, n_t)}
\end{equation}

The trust is computed as the average of all trust traces left within a bin according to Eq. (\ref{eq:bin_trust}), where $trace_{i,j} = \{trace|trace \in Bin_{i,j}\}$ which represents the vehicles traces left within $Bin_{i,j}$ and $\mathbf{T}_{trace_{i,j}}$ represents the trust of those traces.

\begin{equation}\label{eq:bin_trust}
    \hat{T}_{i,j} = Mean(\mathbf{T}_{trace_{i,j}})
\end{equation}

\subsection{Experimental results}

\subsubsection{Static Trust Field}

\begin{figure*}
     \centering
     \begin{subfigure}[b]{0.3\textwidth}
         \centering
         \includegraphics[width=\textwidth]{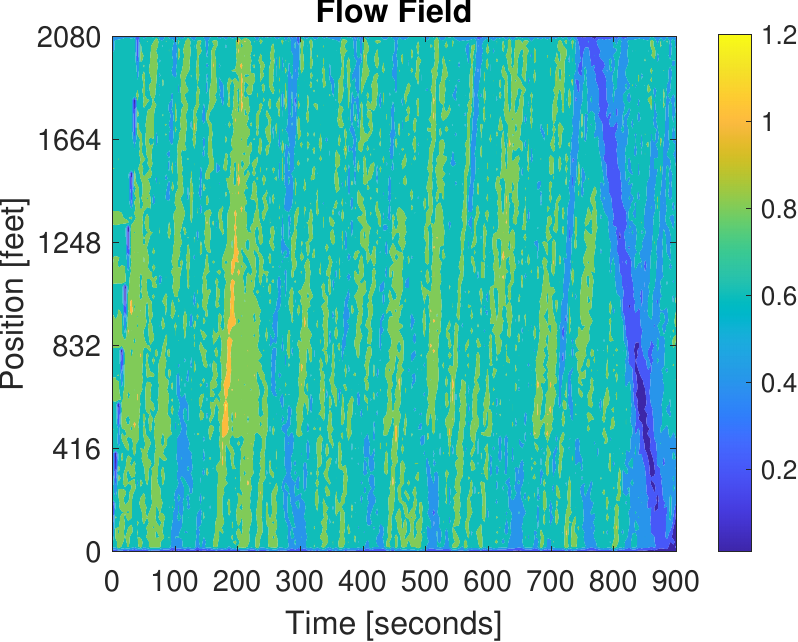}
         \caption{}
         \label{fig:trust_field_ngsim_uni_flow}
     \end{subfigure}
     \hfill
     \begin{subfigure}[b]{0.3\textwidth}
         \centering
         \includegraphics[width=\textwidth]{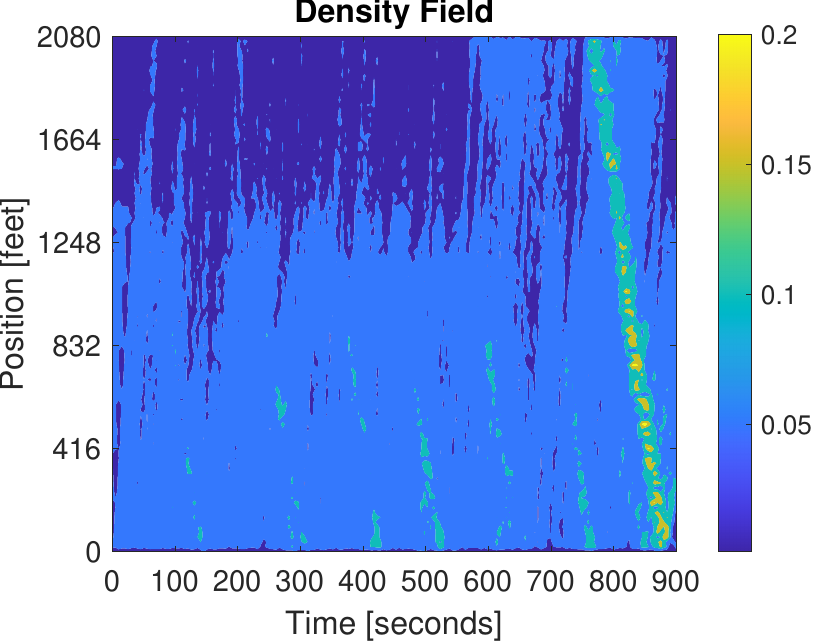}
         \caption{}
         \label{fig:trust_field_ngsim_uni_density}
     \end{subfigure}
     \hfill
     \begin{subfigure}[b]{0.3\textwidth}
         \centering
         \includegraphics[width=\textwidth]{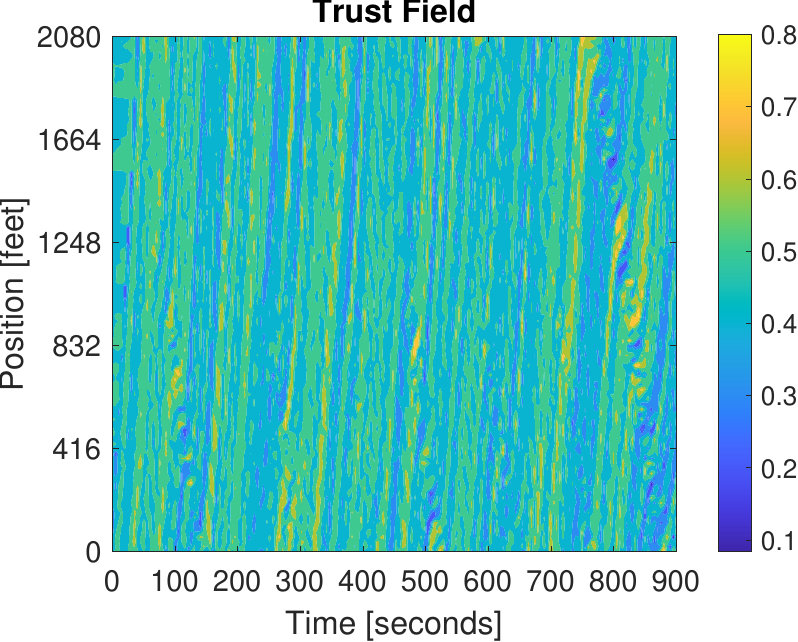}
         \caption{}
         \label{fig:trust_field_ngsim_uni_trust}
     \end{subfigure}
        \caption{Spatiotemporal data for the US-101 highway obtained from binning the NGSIM trajectory data. \textbf{(a)} Spatiotemporal flow data. The flow data is obtained as the product of the velocity and density data sets. \textbf{(b)} Spatiotemporal density data. The periods corresponding to free-flowing traffic have a smaller density, and periods corresponding to traffic jams are a result of high density. \textbf{(c)} Trust field for the NGSIM data. The trust score is constant for each vehicle throughout the simulation and obtained from a uniform distribution.}
        \label{fig:trust_field_ngsim_uni}
\end{figure*}

Fig. \ref{fig:trust_field_ngsim_uni} represents the spatiotemporal data for the NGSIM US-101 freeway obtained from binning the trajectory data. For all the sub-figures, the x-axis represents the duration of the observation, ranging from 0 to 900 seconds, and the y-axis indicates the position along the highway, ranging from 0 to 2080 feet. Fig. \ref{fig:trust_field_ngsim_uni_flow} illustrates the spatiotemporal flow data for the US-101 highway, derived by combining velocity and density data. Yellow and green areas indicate periods and locations of high traffic flow, which suggests free-flowing traffic. Blue regions correspond to low traffic flow, potentially indicating congestion or slow-moving traffic. The flow varies over time and space, showing distinct traffic behavior patterns. Fig. \ref{fig:trust_field_ngsim_uni_density} shows the spatiotemporal density data for the US-101 highway, indicating the concentration of vehicles over time and space. Blue areas indicate low vehicle density, corresponding to free-flowing traffic with fewer vehicles per unit area. Yellow areas show high vehicle density, likely indicating traffic congestion or slow-moving traffic where vehicles are closely packed. The density data reveal specific locations and times prone to high vehicle concentration, highlighting potential congestion zones. Fig. \ref{fig:trust_field_ngsim_uni_trust} represents the trust field obtained from the NGSIM trajectory data using the binning method. In the NGSIM trajectory data, the trust score for each vehicle remained constant throughout the simulation, and the initial trust score is taken from a uniform distribution in the range $[0, 1]$. We observed that the trust field follows a similar pattern to the other two macroscopic traffic data: flow and density fields. The trust field is evolving in time and space and follows the vehicle's trajectory.

\begin{figure*}
     \centering
     \begin{subfigure}[b]{0.3\textwidth}
         \centering
         \includegraphics[width=\textwidth]{Figure_Trust/us101/ngsim_us_101_Flow_field.pdf}
         \caption{}
         \label{fig:trust_field_logit_flow}
     \end{subfigure}
     \hfill
     \begin{subfigure}[b]{0.3\textwidth}
         \centering
         \includegraphics[width=\textwidth]{Figure_Trust/us101/ngsim_us_101_Density_field.pdf}
         \caption{}
         \label{fig:trust_field_logit_density}
     \end{subfigure}
     \hfill
     \begin{subfigure}[b]{0.3\textwidth}
         \centering
         \includegraphics[width=\textwidth]{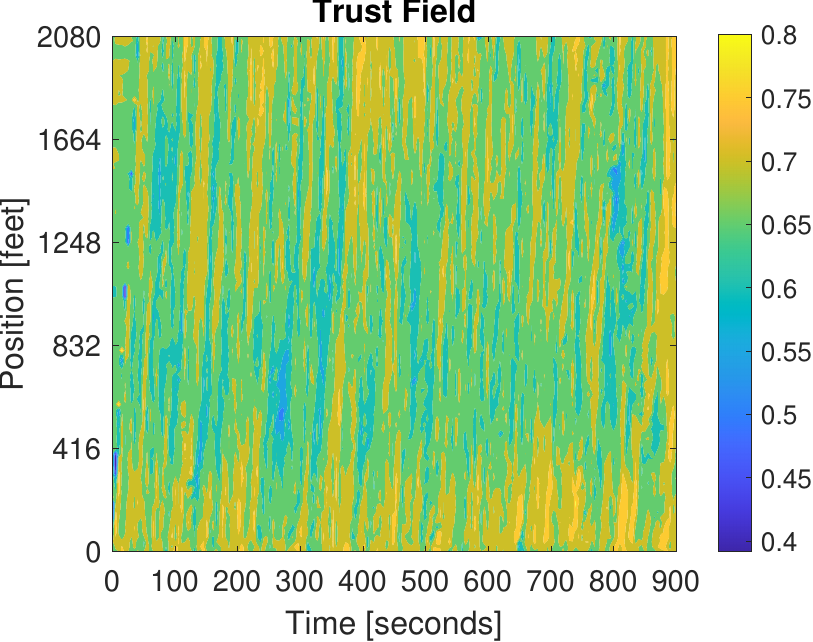}
         \caption{}
         \label{fig:trust_field_logit_trust}
     \end{subfigure}
        \caption{Spatiotemporal data for the US-101 highway obtained from binning the NGSIM trajectory data. \textbf{(a)} Spatiotemporal flow data. The flow data is obtained as the product of the velocity and density data sets. \textbf{(b)} Spatiotemporal density data. The periods corresponding to free-flowing traffic have a smaller density, and periods corresponding to traffic jams are a result of high density. \textbf{(c)} Trust field for the NGSIM data. The trust score evolves dynamically for each vehicle throughout the simulation.}
        \label{fig:trust_field_ngsim_logit}
\end{figure*}

\subsubsection{NGSIM Dynamic Trust}

Fig. \ref{fig:trust_field_logit_flow} and \ref{fig:trust_field_logit_density} shows the same spatiotemporal data as Fig. \ref{fig:trust_field_ngsim_uni_flow} and \ref{fig:trust_field_ngsim_uni_density} for the NGSIM US-101 freeway. This data was obtained by binning the trajectory data. However, the trust field shown in Fig. \ref{fig:trust_field_logit_trust} is obtained from the LogitTrust model. The trust score of each vehicle evolved dynamically and was updated at every second (discussed in Section \ref{sec:trust_evaluation}). Here, $10\%$ of the nodes were malicious, and we can observe that it is reflected in the trust field represented by the blue areas. The yellow and green areas show the trustworthy regions. From the trust field, it is evident that trust is evolving both in space and time.

\section{Conclusion} \label{sec:conc}

This study provides empirical evidence for a trust field in VANETs. It introduces a flexible trust evaluation model that can adjust to the network's operational and environmental fluctuations. Using the LogitTrust model, we have shown that trust fields can adapt and change, accurately representing the reliability of individual nodes by considering their actions and performance over time. The findings suggest that the trust field corresponds with other macroscopic traffic statistics, such as flow and density, and evolves spatiotemporally, thus providing an extensive network reliability assessment. This research contributes to advancing secure and reliable VANET applications, providing a robust mechanism to counteract malicious activities and enhance overall system integrity. In future work, we plan to develop a trust field model mathematically that can evolve both in space and time and integrate a more sophisticated trust computation model to validate and refine the trust field.

\bibliographystyle{IEEEtran}
\bibliography{main}

\end{document}